\documentclass[12pt]{article}
\usepackage{amsmath,amssymb,amsfonts,graphicx,wrapfig}

\newtheorem{theorem}{Theorem}

\newenvironment{proof}{{\sc Proof}}{~\hfill QED}

\newcommand{\li}{{\mathrm{lim\,inf}}}

\begin{document}

{\large
\begin{center}
Existential monadic second order logic of undirected graphs: the Le Bars conjecture is false\footnote{This work is supported by the grant N 16-31-60052 of Russian Foundation for Basic Research and the grant N NSh.-6760.2018.1}
\end{center}
}

\begin{center}
S.N. Popova\footnote{Moscow Institute of Physics and Technology, laboratory of advanced combinatorics and network applications\\ popovaclaire@mail.ru}, M.E. Zhukovskii\footnote{Moscow Institute of Physics and Technology, laboratory of advanced combinatorics and network applications; The Russian Presidential Academy of National Economy and Public Administration\\zhukmax@gmail.com}
\end{center}

\kern 0.5 cm

\begin{center}
Abstract
\end{center}

In 2001, J.-M. Le Bars disproved the zero-one law (that says that every sentence from a certain logic is either true asymptotically almost surely (a.a.s.), or false a.a.s.) for existential monadic second order sentences (EMSO) on undirected graphs. He proved that there exists an EMSO sentence $\phi$ such that ${\sf P}(G_n \models\phi)$ does not converge as $n\to\infty$ (here, the probability distribution is uniform over the set of all graphs on the labeled set of vertices $\{1, \dots, n\}$). In the same paper, he conjectured that, for EMSO sentences with 2 first order variables, the zero-one law holds. In this paper, we disprove this conjecture.

\section{Introduction}
\label{s1}

For undirected graphs, sentences in the first order language (FO sentences) are constructed using relational symbols $\sim$ (interpreted as adjacency) and $=$ (equality of vertices), logical connectives $\neg,\rightarrow,\leftrightarrow,\vee,\wedge$, variables $x,y,z, \ldots$ (with subscripts and superscripts) interpreted as vertices of a graph, quantifiers $\forall,\exists$ and parentheses.  Monadic second order, or MSO, sentences are built of the above symbols of the first order language, as well as the variables $X,Y,Z,\ldots$ (with subscripts and superscripts) that are interpreted as unary predicates. In an MSO sentence, variables $x,y,z, \ldots$ (interpreted as vertices) are called {\it FO variables}, and variables $X,Y,X_1,\ldots$ (that express sets) are called {\it monadic (or MSO) variables}. If, in an MSO sentence $\phi$, all the MSO variables are existential and in the beginning, that is 
\begin{equation}
\phi=\exists X_1\ldots\exists X_m\,\,\varphi(X_1,\ldots,X_m)
\label{general_EMSO}
\end{equation} 
where $\varphi(X_1,\ldots,X_m)$ is a FO sentence with unary predicates $X_1,\ldots,X_m$, then the sentence is called {\it existential} monadic second order (EMSO). Sentences must have finite number of logical connectives. We call the number of nested quantifiers in a longest sequence of nested quantifiers of a formula $\varphi$ {\it the quantifier depth} (note that tautological equivalent formulae may have different quantifier depths). {\it The FO quantifier depth} of an EMSO sentence~(\ref{general_EMSO}) is the quantifier depth of $\varphi(X_1,\ldots,X_m)$. For example, the EMSO sentence
$$
 \exists X \quad[\exists x\exists y\,\,X(x_1)\wedge \neg X(x_2)]\wedge\neg[\exists x \exists y\,\,X(y)\wedge\neg X(z)\wedge y\sim z]
$$
has 2 FO variables, quantifier depth $3$, FO quantifier depth 2 and expresses the property of being disconnected. The quantifier depth of a sentence has the following clear algorithmic interpretation: a FO sentence of quantifier depth $k$ on an $n$-vertex graph can be verified in $O(n^k)$ time. It is very well known (see, e.g.,~\cite{Libkin}, Proposition 6.6) that the same is true for the number of variables: a FO sentence with $k$ variables on an $n$-vertex graph can be verified in $O(n^k)$ time. The later statement is stronger because, clearly, every FO sentence of quantifier depth $k$ may be rewritten using at most $k$ variables.

In what follows, for a sentence $\phi$, we use the usual notation ``$G\models\phi$'' from model theory to denote that $\phi$ is true for $G$.\\

In 1969, Y.V. Glebskii, D.I. Kogan, M.I. Liogon'kii and V.A. Talanov, and independently R. Fagin in 1976~\cite{Fagin}, proved that any FO sentence is either true for almost all graphs, or false for almost all graphs. Clearly, this result can be reformulated in terms of the binomial random graph $G(n,1/2)$. For arbitrary $p$, let $G(n,p) = (V_n,E)$, where $V_n=\{1,\ldots,n\}$, and each pair of vertices is connected by an edge with probability $p$ and independently of other pairs. For more information, we refer readers to the books \cite{AS,Bollobas,Janson}. {\it The zero-one law} of Glebskii et al. and Fagin says that, for every FO sentence $\phi$, either ${\sf P}(G(n,1/2)\models\phi)\to 0$ as $n\to\infty$, or ${\sf P}(G(n,1/2)\models\phi)\to 1$ as $n\to\infty$ (or, in other words, either $\phi$ is true for $G(n,1/2)$ a.a.s., or false a.a.s.). Below, we give a brief history of studying logical laws for this random graph model, for more details (especially, for FO logic) see, e.g.,~\cite{Survey,Strange}. For MSO, the zero-one law for $G(n,1/2)$ was disproved by M. Kaufmann and S. Shelah in 1985~\cite{Kaufmann_Shelah}. They even prove that there is no {\it MSO convergence law} (i.e., there is an MSO sentence $\phi$ such that ${\sf P}(G(n,1/2)\models\phi)$ does not converge). After that, in 1987~\cite{Kaufmann}, Kaufmann proved that there exists an EMSO sentence with 4 binary relations that has no asymptotic probability. The non-convergence result for 1 binary symmetric relation (i.e., for $G(n,1/2)$) was obtained by J.-M. Le Bars in 2001~\cite{Le_Bars}. Note that the construction of Kaufmann has 4 monadic variables and 9 first order variables, and the sentence proposed by Le Bars has even more variables (of both types). In the above mentioned paper, Le Bars conjectured that, for EMSO sentences with 2 FO variables, $G(n,1/2)$ obeys the zero-one law. In this paper, we disprove this conjecture.\\

The paper has the following organization. In Section~\ref{conj}, we construct an EMSO sentence with 2 FO variables and prove that the probability that it is true on $G(n,1/2)$ does not converge. This construction immediately implies that the minimum number of FO variables of an EMSO sentence without convergence equals 2, and the same is true for the FO quantifier depth. In Section~\ref{related}, we prove that our construction is, in some sense, best possible. Namely, we rewrite the sentence in certain ways that exploit only 1 monadic variable. These tautological equivalents of the sentence still have small number of FO variables. Upon proving this, we conclude that even 1 monadic variable and 2 FO variables (in one sentence) are enough for non-convergence, but the minimum FO quantifier depth of an EMSO sentence having 1 monadic variables and without convergence equals 3. Moreover, in the same section, we prove that there is a dense subset $\mathcal{P}\subset[0,1]$ such that, for every $p\in\mathcal{P}$, there is an EMSO sentence $\phi$ with 2 FO variables such that ${\sf P}(G(n,p)\models\phi)$ does not converge as $n\to\infty$.

\section{Disproving the conjecture}
\label{conj}

This section is devoted to the proof of the following result.

\begin{theorem}
There exists an EMSO sentence $\phi$ with two monadic variables and FO quantifier depth 2 such that the probability ${\sf P}(G(n, 1/2) \models\phi)$ does not converge as $n \to \infty$.
\label{one}
\end{theorem}

\proof

Let
\begin{multline*}
\phi=\exists X_1 \exists X_2 \quad 
[\exists x\,\, X_1(x)]\land[\exists x\,\, X_2(x)]\land\\
\varnothing(X_1,X_2)\land 
\mathrm{Cl}(X_1) \land 
\mathrm{Cl}(X_2)  \land 
\mathrm{lnd}(X_1, X_2) \land 
\mathrm{N}(X_1 \vee X_2) \land
\mathrm{Max}(X_1,X_2),
\end{multline*}
where
$$
 \varnothing(X_1,X_2)=\neg[\exists x\,\,X_1(x) \wedge X_2(x)],
$$
$$
\mathrm{Cl}(X) = \forall x \forall y\quad [X(x) \land X(y) \land (x \neq y)] \to [x \sim y],
$$
$$
 \mathrm{Ind}(X, Y) = \forall x \forall y\quad [X(x) \land Y(y)]\to [x \not \sim y],
$$
$$
\mathrm{N}(X) = \exists y \quad [\neg X(y)] \land [\forall x\,\, (X(x) \to  y \sim x)],
$$
$$
\mathrm{Max}(X_1, X_2) = \forall y\quad \neg[X_1(y)\vee X_2(y)]\to [(\exists x \,\,X_1(x) \land x\sim y)\land(\exists x \,\,X_2(x) \land x\sim y)].
$$
In other words, $\phi$ says that there two disjoint cliques such that
\begin{itemize}
\item there are no edges between them,
\item there is a common neighbor of vertices of both cliques,
\item every vertex outside both cliques has neighbors in both.
\end{itemize}

Let $X(k, l)$ be the number of triples $(X_1, X_2, x)$ such that
\begin{itemize}
\item $X_1$ and $X_2$ are two non-empty cliques of sizes $k,l$ respectively in $G(n,\frac{1}{2})$, 
\item $X_1 \cap X_2 = \varnothing$,
\item $x \notin X_1 \cup X_2$,
\item there are no edges between $X_1$ and $X_2$, 
\item $x$ is adjacent to every vertex of $X_1 \cup X_2$,
\item every vertex from  $\overline{X_1 \cup X_2}$ has a neighbor in $X_1$ and a neighbor in $X_2$.
\end{itemize}
It is immediate that ${\sf P}(G(n,1/2)\models\phi) = {\sf P}(\exists k \exists l: \quad X(k, l) > 0)$.  Let us estimate ${\sf E}X(k, l)$.
\begin{multline*}
{\sf E}X(k, l) = \frac{n!}{k! l! (n - k - l)!} (n - k - l) 2^{-\frac{(k + l)^2}{2} + \frac{k + l}{2}} 2^{-(k + l)} \cdot 
 (1 - 2^{-k})^{n - k - l - 1} (1 - 2^{-l})^{n - k - l - 1}.
\end{multline*}
Thus,
\begin{equation}
{\sf E}X(k, l) \le \frac{n^{k + l + 1} e^{k + l}}{k^k l^l} 2^{-\frac{(k + l)^2}{2} - \frac{k + l}{2}} e^{-(n - k - l - 1)(2^{-k} + 2^{-l})}.
\label{expectation}
\end{equation}
Set
$$
f(k, l) = (k + l + 1) \ln n - k \ln k - l \ln l + k + l - \frac{(k + l)^2}{2} \ln 2 - \frac{k + l}{2} \ln 2 - (n - k - l) (2^{-k} + 2^{-l}).
$$ 
Let us prove that $f$ has a maximum and estimate it. Compute the derivatives of $f$:
\begin{equation}
\frac{\partial f}{\partial k}(k, l) = \ln n - \ln k - (k+l) \ln 2 - \frac{1}{2} \ln 2 + (n - k - l) 2^{-k} \ln 2 + 2^{-k} + 2^{-l},
\label{partial1}
\end{equation}
\begin{equation}
\frac{\partial f}{\partial l}(k, l) = \ln n - \ln l - (k+l) \ln 2 - \frac{1}{2} \ln 2 + (n - k - l) 2^{-l} \ln 2 + 2^{-k} + 2^{-l},
\label{partial2}
\end{equation}
$$
\frac{\partial^2 f}{\partial k^2}(k, l) = -\frac{1}{k} - \ln 2 - (n - k - l) 2^{-k} \ln^2 2 - 2^{1-k} \ln 2,
$$
$$
\frac{\partial^2 f}{\partial l^2}(k, l) = -\frac{1}{l} - \ln 2 - (n - k - l) 2^{-l} \ln^2 2 - 2^{1-l} \ln 2,
$$
$$
\frac{\partial^2 f}{\partial k \partial l}(k, l) = -\ln 2 - (2^{-k} + 2^{-l}) \ln 2.
$$

The matrix $A$ of second-order partial derivatives of $f$ is negative definite for all real $k,l\geq 1$ such that $k+l\leq n-1$ (since $\frac{\partial^2 f}{\partial k^2}<0$ and the determinant $|A|=\frac{\partial^2 f}{\partial k^2}\frac{\partial^2 f}{\partial l^2}-\left(\frac{\partial^2 f}{\partial k\partial l}\right)^2>0$). Let us find (the only) zero of $f'=(\frac{\partial f}{\partial k},\frac{\partial f}{\partial l})$. Due to symmetry reasons, it should be of the form $(k^*,k^*)$. From~(\ref{partial1}),~(\ref{partial2}), it is easy to see that $k^*=\frac{\ln n - \ln \ln n + \ln\ln 2}{\ln 2}+O(\frac{\ln \ln n}{\ln n})$. Then
\begin{multline*}
f(k^*, k^*) = \frac{2}{\ln 2} (\ln n - \ln \ln n + \ln \ln 2) \ln n + \ln n - \frac{2}{\ln 2} (\ln n - \ln \ln n) (\ln \ln n - \ln \ln 2) + \\ +
\frac{2}{\ln 2} \ln n
- \frac{2}{\ln 2} (\ln n - \ln \ln n + \ln \ln 2)^2 - \ln n - \frac{2}{\ln 2} \ln n + O(\ln \ln n) = O(\ln \ln n).
\end{multline*}

Now, let us find a sequence $n_k$ such that ${\sf P}(G(n_k,1/2)\models\phi)\to 0$ as $k\to\infty$.

For $k\in\mathbb{N}$, set $n=n_k=\lfloor 2^{k + \frac{1}{2}} k \rfloor$. Then $k^* =k^*(n)= k + \frac{1}{2} + o(1)$. Let $\tilde k, \tilde l \in \mathbb N$. Then $\tilde k = k^* + \Delta_k$, $\tilde l = k^* + \Delta_l$, where
$|\Delta_k| \ge \frac{1}{2} + o(1)$, $|\Delta_l| \ge \frac{1}{2} + o(1)$. We have
\begin{multline*}
f(\tilde k, \tilde l) - f(k^*, k^*) \le -\left((\Delta_k + \Delta_l) \ln n + \frac{\ln n}{\ln 2} (2^{-\Delta_k} + 2^{-\Delta_l} - 2)\right) (1 + o(1)) = \\
= -\frac{\ln n}{\ln 2} (g(\Delta_k \ln 2) + g(\Delta_l \ln 2)) (1 + o(1)), 
\end{multline*}
where
$g(x) = x + e^{-x} - 1$, and the bound is uniform over all $\tilde k,\tilde l$. Therefore, $f(\tilde k, \tilde l)-f(k^*,k^*) \le - c \ln n$ for any $\tilde k, \tilde l \in \mathbb N$ and some constant $c > 0$.
Moreover, for all $n$ large enough, we have $f(\tilde k, \tilde l)-f(k^*,k^*) \le -3\ln n$ for any $\tilde k, \tilde l \in \mathbb N$ such that either $|\tilde k - k^*|\ge 5$ or $|\tilde l - k^*| \ge 5$.
Hence, from~(\ref{expectation}), for $n$ large enough,
$$
{\sf P}(\exists \tilde k \exists \tilde l: \quad X(\tilde k, \tilde l) > 0) \le\sum_{\tilde k,\tilde l}{\sf E}X(\tilde k,\tilde l)\leq\sum_{\tilde k,\tilde l}e^{f(\tilde k,\tilde l)+1}\leq 
$$
$$
e^{-3\ln n+1+f(k^*,k^*)} n^2 + 100 e^{-c\ln n+1+f(k^*,k^*)}  = o(1).
$$
Therefore, we have 
\begin{equation}
\lim_{k \to \infty} {\sf P}\left(G\left(n_k,\frac{1}{2}\right)\models\phi\right) = 0.
\label{zero_limit}
\end{equation}

Finally, let us find a sequence $n_k$ such that $\lim\inf_{k\to\infty}{\sf P}(G(n_k,1/2)\models\phi)>0$.

Set $n = n_k = 2^k k$. Let us estimate ${\sf E}X(k, k)$:
$$
{\sf E}X(k, k) = \frac{n^{2k+1} e^{2k}}{k^{2k} 2\pi k} 2^{-2k^2 - k} e^{-2n 2^{-k}} (1 + o(1)) = \frac{1}{2\pi}(1 + o(1)).
$$
It can be shown similarly as in the above proof that
$$
{\sf P}\biggl(\exists \tilde k \exists \tilde l  \quad  [(\tilde k, \tilde l) \neq (k, k)] \land [X(\tilde k, \tilde l) > 0]\biggr)  = o(1).
$$
Let us estimate ${\sf P}(X(k, k) > 0)$. Consider the random variable $\tilde X(k, k) = X(k, k)/2$ which counts the number of triplets $(X_1,X_2,x)$ where the order of cliques $X_1,X_2$ in the triplet does not matter.

Let us remind that
$$
 {\sf E}\tilde X(k,k)=\frac{1}{2}{n\choose 2k}{2k\choose k}(n-2k)2^{-{2k\choose 2}}2^{-2k} e^{-2n 2^{-k}+o(1)}.
$$
Below, we compute ${\sf E} \tilde X^2(k, k)$ in the usual way: it equals to the summation of probabilities of the events that triplets $(X_1^1,X_2^1,x^1)$ and $(X_1^2,X_2^2,x^2)$ consisting of two disjoint $k$-sets $X_1^j,X_2^j$ and a vertex $x^j$ outside them have the above property (both sets induce cliques, there are no edges between them, the distinguished vertex is a common neighbor of them, and every other vertex have neighbors in every set) over all possible pairs of such triplets.

We distinguish four cases w.r.t. possible intersections of elements of two triplets.

The contribution of pairs of equal triplets is ${\sf E}\tilde X(k,k)$. 

The contribution of pairs of triplets sharing both sets (but not a vertex) equals 
$$
\frac{n!}{2 k! k! (n - 2k)!} (n - 2k) (n - 2k - 1) 2^{- {{2k} \choose 2}} 2^{-4k} e^{-2n 2^{-k}} (1 + o(1))=o(1).
$$

The contribution of pairs of triplets with disjoint sets (i.e., $X_{i_1}^{j_1}\cap X_{i_2}^{j_2}=\varnothing$ whenever  $(i_1,j_1)\neq(i_2,j_2$)) equals
$$
 \frac{n!}{4 (k!)^4 (n - 4k)!} n^2 2^{- 2 {{2k} \choose 2}} 2^{-4k} e^{-4n 2^{-k}} (1 + o(1))= ({\sf E} \tilde X(k, k))^2(1+o(1)).
$$

Finally, let us denote the contribution of triplets such that $|(X_1^1\sqcup X_2^1)\cap(X_1^2\sqcup X_2^2)|=j$ by $A_j$. Then, clearly,
$$
{\sf E} \tilde X^2(k, k) = {\sf E} \tilde X(k, k) + ({\sf E} \tilde X(k, k))^2 + \sum_{j = 1}^{2k - 1} A_j + o(1).
$$
For estimating $A_j$, $j\in\{1,\ldots,2k-1\}$, let us consider two sets of $2k$ vertices $U^1$ and $U^2$ in $V_n$ that share $j$ vertices. Let us assume that each set $U^i$ has two parts $U^i_1$ and $U^i_2$ of the same size $k$, and, for $i\in\{1,2\}$, the $i$-th part $U^1_i$ of the first set has $j_i\leq k$ common vertices with the $i$-th part $U^2_i$ of the second set. Moreover, let $U^1_1\cap U^2_2=U^1_2\cap U^2_1=\varnothing$. So, $j_1+j_2=j$. Then the probability that a vertex outside $U^1\cup U^2$ has neighbors in each of the sets $U^1_1$, $U^1_2$, $U^2_1$ and $U^2_2$ equals
$$
 \left(1-2^{-j_1}\right)\left(1-2^{-j_2}\right)+2^{-j_1}\left(1-2^{-j_2}\right)\left(1-2^{-(k-j_1)}\right)^2+2^{-j_2}\left(1-2^{-j_1}\right)\left(1-2^{-(k-j_2)}\right)^2+
$$ 
\begin{equation} 
 2^{-j}\left(1-2^{-(k-j_1)}\right)^2\left(1-2^{-(k-j_2)}\right)^2=1+2^{-k}\left(-4+2^{-k+j_1}+2^{-k+j_2}+O(2^{-k})\right).
\label{aux_probab}
\end{equation}
As $2^{-k+j_1}+2^{-k+j_2}$ is convex in $j_2$ (recall that $j_1=j-j_2$ and $j$ is fixed), then it achieves its maximum in one of the endpoints of the interval of the admissible values of $j_2$ (no matter which one, since the expression is symmetric w.r.t. $j_1,j_2$). In this way, we get that the expression to the right in~(\ref{aux_probab}) is at most $1+2^{-k}(-4+2^{-k+j}+O(2^{-k}))$ if $j\leq k$ and at most $1+2^{-k}(-3+2^{-2k+j}+O(2^{-k}))$ if $j>k$.

Then, the probability that every vertex outside $U^1\cup U^2$ has neighbors in each of the sets $U^1_1$, $U^1_2$, $U^2_1$ and $U^2_2$ is at most
$$
 \left(1+2^{-k}(-4+2^{-k+j}+O(2^{-k}))\right)^{n-4k+j}\leq e^{n2^{-k}(-4+2^{-k+j}+O(2^{-k}))+o(1)}=e^{-4k+k2^{-k+j}+o(1)}
$$
if $j\leq k$, and at most
$$
 \left(1+2^{-k}(-3+2^{-2k+j}+O(2^{-k}))\right)^{n-4k+j}\leq e^{n2^{-k}(-3+2^{-2k+j}+O(2^{-k}))+o(1)}=e^{-3k+k2^{-2k+j}+o(1)}
$$
if $j>k$.

Let $B_j$ be the expected number of pairs of vertices (not necessarily distinct) $(u_1,u_2)$ such that $u_1$ has neighbors both in $U_1^1$ and $U_2^1$, $u_2$ has neighbors both in $U_1^2$ and $U_2^2$. Clearly,
$$
B_j \le n^2 2^{-4k} + n 2^{-(4k - j)} + 2n (2k - j) 2^{-(4k - j)} + (2k - j)^2 2^{-2(2k - j)}.
$$

From the above estimations and the fact the number of ways of dividing a set of $2k-j$ vertices into two parts of given sizes is at most  ${{2k - j} \choose {k - \lfloor j/2 \rfloor}}$, we get
$$
A_j \le {n \choose {2k}} {{2k} \choose k} {2k \choose j} {{n - 2k} \choose 2k-j} 2^{-2 {{2k} \choose 2} + {j \choose 2}} B_j{{2k - j} \choose {k - \lfloor j/2 \rfloor}}  e^{-4k} e^{k C_j} (1 + o(1)),
$$
where $C_j = 2^{-(k - j)}$ for any $j \le k$ and $C_j = 1 + 2^{-(2k - j)}$ for any $j > k$. Therefore, 
$$
A_j \le 4({\sf E} \tilde X(k, k))^2 \cdot \frac{{2k \choose j} {{n - 2k} \choose 2k-j}
{{2k - j} \choose {k - \lfloor j/2 \rfloor}}}{{n \choose 2k} {2k \choose k} 
n^2 2^{-4k}} 2^{{j \choose 2}} e^{k C_j} B_j (1+o(1)).
$$
Set
$$
F_j = {2k \choose j} {{n - 2k} \choose 2k-j} {{2k - j} \choose {k - \lfloor j/2 \rfloor}} 2^{{j \choose 2}}.
$$
Then 
$$
 \frac{A_j}{({\sf E} \tilde X(k, k))^2}\le \frac{4F_jB_je^{kC_j}}{{n \choose 2k} {2k \choose k} n^2 2^{-4k}}(1+o(1)).
$$
It remains to bound from above the right side of the last inequality. We will do that separately for $j\in\{1,\ldots,4\}$, $j\in\{5,\ldots,2k-5\}$ and $j\in\{2k-4,\ldots,2k-1\}$. 

Let $j\in\{1,\ldots,5\}$. Then
$$
F_j = n^{2k - j} (2k)^j \frac{2^{{j \choose 2}} {{2k - j} \choose {k - \lfloor j/2 \rfloor}}}{(2k-j)!j!} (1 + o(1)),\quad B_j\leq n^2 2^{-4k}(1+o(1)),\quad e^{kC_j}=1+o(1).
$$
In this case,
\begin{equation}
 \frac{4F_jB_je^{kC_j}}{{n \choose 2k} {2k \choose k} n^2 2^{-4k}}=O\left(\frac{n^{2k - j} (2k)^j{{2k -j} \choose {k - \lfloor j/2 \rfloor}}}{{n \choose 2k} {2k \choose k}(2k-j)!}\right)=O\left(\frac{n^{2k-j}(2k)^{2k+j}e^{2k-j}}{n^{2k}e^{2k}(2k-j)^{2k-j}}\right)=
 O\left(\frac{k^{2j}}{n^j}\right).
\label{small_j}
\end{equation}
Moreover,
$$
F_{2k-j} = n^j \frac{(2k)^j {j \choose \lceil{ j/2 \rceil} }}{(j!)^2} e^{k (1 + 2^{-j})} 2^{2k^2-k(2j+1)+\frac{j(j+1)}{2}} (1 + o(1)),\quad B_{2k-j}\leq 1,\quad e^{kC_{2k-j}}=e^{k(1+2^{-j})};
$$
$$
 \frac{4F_{2k-j}B_{2k-j}e^{kC_{2k-j}}}{{n \choose 2k} {2k \choose k} n^2 2^{-4k}}=O\left(\frac{n^j (2k)^j e^{k (1 + 2^{-j})} 2^{2k^2-k(2j+1)} }{{n \choose 2k} {2k \choose k}n^2 2^{-4k}}\right)=O\left(\frac{(k2^k)^{j}(2k)^{2k+j}e^{k(1+2^{-j})}2^{2k^2}}{(k2^k)^{2k}e^{2k}k^{3/2}2^{k(2j+1)}}\right)=
$$ 
\begin{equation}
 O\left(\frac{k^{2j-3/2}}{e^{k(1-2^{-j})}2^{k(j-1)}}\right).
\label{large_j}
\end{equation} 

It immediately follows that
\begin{equation}
\sum_{j\in\{1,\ldots,4\}\sqcup\{2k-4,\ldots,2k-1\}}\frac{A_j}{({\sf E} \tilde X(k, k))^2}=o(1).
\label{close_to_the_bounds}
\end{equation}

Since $B_j$ and $C_j$ increase as $j$ grows from $1$ to $2k-1$, the following bounds hold for all $j$: $B_j\leq 1$, $e^{kC_j}\leq e^{3k/2}$. Consequently, for $j\in\{5,\ldots,2k-5\}$, we may use the following bound:
$$
 \frac{A_j}{({\sf E} \tilde X(k, k))^2}\le \frac{4F_j e^{3k/2}}{{n \choose 2k} {2k \choose k} n^2 2^{-4k}}(1+o(1)).
$$
It is very easy to check that $F_j$ first decreases and then increases as $j$ grows from $1$ to $2k-1$. Therefore, for every $j\in\{5,\ldots,2k-5\}$, the estimates~(\ref{small_j})~and~(\ref{large_j}) imply
$$
 \frac{A_j}{({\sf E} \tilde X(k, k))^2}\leq\max\{F_5,F_{2k-5}\}\frac{4e^{3k/2}}{{n \choose 2k} {2k \choose k} n^2 2^{-4k}}(1+o(1))=
$$
$$ 
 \max\left\{\frac{e^{3k/2}}{n^2 2^{-4k}}O\left(\frac{k^{10}}{n^5}\right),e^{k(1-2^{-5})}O\left(\frac{k^{17/2}}{e^{k(1-2^{-5})}16^{k}}\right)\right\}=O\left(\frac{k^3}{(8/e^{3/2})^k}\right),
$$
and this bound is uniform. From this bound and~(\ref{close_to_the_bounds}), we get $\sum_{j = 1}^{2k - 1} A_j = o(({\sf E} \tilde X(k, k))^2)$ and ${\sf E} \tilde X^2(k, k) =  {\sf E} \tilde X(k, k) +  ({\sf E} \tilde X(k, k))^2(1+o(1)) + o(1)$.

Finally,
$$
{\sf P}(\tilde X(k, k) > 0) \ge \frac{({\sf E} \tilde X(k, k))^2}{{\sf E} \tilde X^2(k, k)} = \frac{1}{4\pi + 1}(1 + o(1)).
$$
Therefore, by the Paley--Zygmund inequality~\cite{PZ1,PZ2},
\begin{equation}
\li_{k \to \infty} {\sf P}\left(G\left(n_k,\frac{1}{2}\right)\models\phi\right)>0.
\label{positive_limit}
\end{equation}

From~(\ref{zero_limit})~and~(\ref{positive_limit}), the ${\sf P}(G(n,1/2)\models\phi)$ does not converge as $n\to\infty$. $\Box$\\

\section{Related results}
\label{related}

\subsection{Minimum quantifier depth and minimum number of variables}

It is not difficult to show using the modification of Ehrenfeucht-Fra\"{\i}ss\'e games for EMSO (see, e.g.,~\cite{Zhuk_Le_Bars}) which is also known as the Fagin game (see~\cite{Libkin}, Chapter 7.3), that FO quantifier depth 1 as well as 1 FO variable are not enough for non-convergence (in the latter case, the variant of the game with one pebble should be considered,~\cite{Libkin}, Chapter 11.4; however, the winning strategy of Duplicator is still obvious in this case). In fact, in both cases, 0-1 law holds. Indeed, if two graphs in the game are large enough (both have more than $2^k$ vertices), then Duplicator wins the game with $k$ set moves and $1$ vertex move (or $1$ pebble).

Therefore, Theorem~\ref{one} implies that the minimum number of FO variables of an EMSO sentence without convergence equals 2. The same is true for the FO quantifier depth. But is the same true if we restrict ourselves to 1 monadic variable? 

Below, we write two tautological equivalents of $\phi$ both with one monadic variable. The first sentence $\phi_1$ has 2 FO variables, and the second sentence $\phi_2$ has the FO quantifier depth 3:
$$
 \phi_1=\exists X\quad(\forall x\forall y \,\, [X(x)\wedge X(y)\wedge(x\neq y)]\to[x\sim y])\wedge
$$
$$
 (\forall x\,\,[\forall y\,\, X(y)\to(x\nsim y)]\to[\forall y\,\,(\forall x\,\, X(x)\to (x\nsim y))\to(x=y\vee x\sim y)])\wedge 
$$
$$
 (\forall x\,\,[\forall y\,\,(\forall x\,\, X(x)\to x\nsim y)\to(x\nsim y)]\to X(x))\wedge
$$
$$
 (\exists x\,\,[\forall y\,\, X(y)\to x\sim y]\wedge[\forall y\,\,(\forall x\,\, X(x)\to x\nsim y)\to x\sim y]).
$$
In order to clarify the construction, let us note that the major part of the sentence $\phi_1$ (excluding quantification over the only monadic variable) is the conjunction of 4 FO formulae. The first one says that the set $[X]$ of vertices $x$ satisfying $X$ induces a clique. The second formula says that the set of vertices (say, $Y$) having no neighbors in $[X]$ induces a clique. The third formula says that every vertex having no neighbors in $Y$ belongs to $[X]$. Finally, the last formula says that there exists a common neighbor of all the vertices in $[X]\sqcup Y$.
$$
 \phi_2=\exists X\quad(\forall x\forall y \,\, [X(x)\wedge X(y)\wedge(x\neq y)]\to[x\sim y])\wedge
$$
$$
 (\forall x\,\,[\forall y\,\, X(y)\to(x\nsim y)]\to[\forall y\,\,(\forall z\,\, X(z)\to (z\nsim y))\to(x=y\vee x\sim y)])\wedge 
$$
$$
 (\forall x\,\,[\forall y\,\,(\forall z\,\, X(z)\to z\nsim y)\to(x\nsim y)]\to X(x))\wedge
$$
$$
 (\exists x\,\,[\forall y\,\, X(y)\to x\sim y]\wedge[\forall y\,\,(\forall z\,\, X(z)\to z\nsim y)\to x\sim y]).
$$
Every formula from the conjunction in $\phi_2$ is tautologically equivalent to the respective formula in $\phi_1$.

In~\cite{Zhuk_Le_Bars}, it is proven that, for EMSO sentences with 1 monadic variable and FO quantifier depth 2, $G(n,1/2)$ obeys 0-1 law. From this, clearly, the minimum FO quantifier depth of an EMSO sentence with 1 monadic variable and without convergence equals 3, and the minimum number of FO variables equals 2.

\subsection{Other $p$}

The proof of Theorem~\ref{one} works only in the case $p=\frac{1}{2}$. In this section, we try to prove non-convergence for other constant $p$. We can not prove it for all $p\in(0,1)$. However, we obtain the following result.

\begin{theorem}
There exists a dense subset $\mathcal P \subset (0, 1)$ such that, for every $p \in \mathcal P$, there exists an EMSO sentence $\phi$ with FO quantifier depth 2 such that the probability ${\sf P}(G(n, p) \models \phi)$ does not converge as $n \to \infty$.
\end{theorem}

{\it Proof}. Surely, if we get non-convergence for certain $p$, then we immediately have the same for $1-p$. To see this it is enough to add a negation in front of every adjacency relation in $\phi$ for which ${\sf P}(G(n, p) \models \phi)$ does not converge. Therefore, we may restrict ourselves with $p\in(0,1/2)$. 

Fix a parameter $\gamma>1$ and consider the equation $p=(1-p)^{\gamma}$. It has the only root in $(0,1)$, and this root is less than $1/2$. 
By the implicit function theorem, the root $p=p(\gamma)$ is a continuous function on $(1,\infty)$. Moreover, it decreases to $0$ as $\gamma$ increases to infinity. Therefore, it is enough to find a dense subset $\Gamma\subset(1,\infty)$ such that, for every $\gamma\in\Gamma$, $G(n,p(\gamma))$ does not obey EMSO convergence law. 

To do this, consider a rational number $\frac{u}{v}>1$, where $\frac{u}{v}$ is an irreducible fraction. Fix a positive integer $h$. Denote $a=2uh$, $t=h(2u-v)$. Consider an EMSO sentence $\phi_h$ with 2 FO variables that expresses the property of existence of $a$ non-trivial disjoint cliques $X_1,\ldots,X_{a}$ with the following properties: 1) there are no edges between them, 2) there is a vertex $x$ which is a common neighbor of vertices of exactly $a/2$ cliques, 3) there is no $i\in\{1,\ldots,a\}$ such that every vertex outside $V(X_1)\cup\ldots\cup V(X_a)\cup\{x\}$ is adjacent to every vertex of the cliques $X_{1+(i\mod a)},\ldots,X_{1+(i+a-t-1\mod a)}$, and has non-neighbors among vertices of every other clique. 

We claim that, for every $\varepsilon>0$, there exist $h\in\mathbb{N}$ and a rational number $\gamma>1$ such that $|\gamma-\frac{u}{v}|<\varepsilon$, and ${\sf P}(G(n,p(\gamma))\models\phi_h)$ does not converge. The desired $\gamma$ is equal to $\frac{a-1}{2a-2t-1}$ (this value approaches $\frac{u}{v}$ as $h\to\infty$).  The rest of the proof is very similar to the proof of Theorem~\ref{one}. Therefore, we only sketch it.

Let $X(k_1,\ldots,k_a)$ be the number of tuples $(X_1,\ldots,X_a, x)$ such that
\begin{itemize}
\item $X_1,\ldots,X_a$ are non-empty cliques of sizes $k_1,\ldots,k_a$ respectively in $G(n,\frac{1}{2})$, 
\item for $i\neq j$, $X_i \cap X_j = \varnothing$,
\item $x \notin X_1 \cup\ldots\cup X_a$,
\item for $i\neq j$, there are no edges between $X_i$ and $X_j$, 
\item there are exactly $a/2$ cliques among $X_1,\ldots,X_a$ such that $x$ is adjacent to every vertex of them,
\item the above property 3) holds.
\end{itemize}
Then ${\sf P}(G(n,p)\models\phi) = {\sf P}(\exists k_1\ldots \exists k_a \quad X(k_1,\ldots, k_a) > 0)$. In a similar way, as in the proof of Theorem~\ref{one}, one can show that ${\sf E}X(k_1,\ldots,k_a)=e^{f(k_1,\ldots,k_a)}$, where $f$ achieves the global maximum in $(k^*,\ldots,k^*)$, where $k^*=\frac{\ln n-\ln\ln n+\ln((a-t)\ln(1/p))}{(a-t)\ln(1/p)}$. Moreover, $f(k^*,\ldots,k^*)=\left(1-\frac{a\ln(1/(1-p))}{2(a-t)\ln(1/p)}\right)\ln n+O(\ln^2\ln n)$. Note that, first, the coefficient $1-\frac{a\ln(1/(1-p))}{2(a-t)\ln(1/p)}$ is positive and approaches $0$ as $h\to\infty$. Second, for every $c>0$, there exists $\delta=\delta(c)$ such that $|k_1-k^*|+\ldots+|k_a-k^*|\geq c$ implies $f(k_1,\ldots,k_a)\leq f(k^*,\ldots,k^*)-\delta\ln n$, and this $\delta(c)$ approaches $\infty$ as $c\to\infty$. Therefore, for $h$ large enough, 
$$
{\sf P}\left(\exists (k_1\ldots k_a)\quad \left[|k_1-k^*|+\ldots+|k_a-k^*|\geq \frac{1}{3}\right] \wedge [X(k_1,\ldots, k_a) > 0]\right)=0.
$$
It remains to consider two sequences $n_1(k)=\lfloor k(1/p)^{(a-t)(k+1/2)}\rfloor$ and $n_2(k)=\lfloor k(1/p)^{(a-t)k}\rfloor$. For the first sequence, $k^*=k+\frac{1}{2}+o(1)$, and therefore, ${\sf P}(G(n_1(k),p)\models\phi)\to 0$ as $k\to\infty$. For the second sequence, $k^*=k+o(1)$, and so, ${\sf E}X(k,\ldots,k)\to\infty$. In the same way as in the proof of Theorem~\ref{one}, it can be shown that ${\sf D}X(k,\ldots,k)=o(({\sf E}X(k,\ldots,k))^2)$. From Chebyshev's inequality, we immediately get that ${\sf P}(G(n_2(k),p)\models\phi)\to 1$ as $k\to\infty$. $\Box$

\renewcommand{\refname}{References}


\begin{thebibliography}{50}
\bibitem{AS} N. Alon, J.H. Spencer, {\it The Probabilistic Method}, John Wiley \& Sons, 2000.
\bibitem{Le_Bars} J.-M. Le Bars, {\it The 0-1 law fails for monadic existential second-order logic on undirected graphs},  Information Processing Letters {\bf 77} (2001) 43--48.
\bibitem{Bollobas} B. Bollob\'{a}s, {\it Random Graphs}, 2nd Edition, Cambridge University Press, 2001.
\bibitem{Fagin} R. Fagin, \emph{Probabilities in finite models}, J. Symbolic Logic {\bf 41} (1976): 50--58.
\bibitem{Glebskii} Y.V. Glebskii, D.I. Kogan, M.I. Liogon'kii, V.A. Talanov, {\it Range and degree of realizability of formulas the restricted predicate calculus}, Cybernetics {\bf 5} (1969) 142--154 (Russian original: Kibernetica {\bf 2}, 17--27).
\bibitem{Janson} S. Janson, T. Luczak, A. Rucinski, {\it Random Graphs}, New York, Wiley, 2000.
\bibitem{Kaufmann} M. Kaufmann, {\it Counterexample to the 0–1 law for existential monadic second-order logic}, Technical Report, CLI Internal Note 32, Computational Logic Inc., December 1987.
\bibitem{Kaufmann_Shelah} M. Kaufmann, S. Shelah, {\it On random models of finite power and monadic logic}, Discrete Math. {\bf 54} (1985) 285--293.
\bibitem{Libkin} L.~Libkin, {\it Elements of finite model theory},  Texts in Theoretical Computer Science. An EATCS Series, Springer-Verlag Berlin Heidelberg, 2004.
\bibitem{PZ1} R. Paley, A. Zygmund, {\it A note on analytic functions in the unit circle}, In Mathematical Proceedings of the Cambridge Philosophical Society, Vol. 28 (1932), P. 266--272.
\bibitem{PZ2} R. Paley, A. Zygmund, {\it On some series of functions, (3)}, In Mathematical Proceedings of the Cambridge Philosophical Society. Vol. 28 (1932), P. 190--205.
\bibitem{Survey} A.M.~Raigorodskii, M.E.~Zhukovskii, {\it Random graphs: models and asymptotic characteristics}, Russian Mathematical Surveys {\bf 70}:1 (2015)~33--81.
\bibitem{Shelah} S. Shelah, J.H. Spencer, \emph{Zero-one laws for sparse random graphs}, J. Amer. Math. Soc., 1988, 1: 97-115.
\bibitem{Strange} J.H.~Spencer, {\it The Strange Logic of Random Graphs}, Springer Verlag, 2001.
\bibitem{Tyszk} J.~Tyszkiewicz, {\it On Asymptotic Probabilities of Monadic Second Order Properties}, Lecture Notes in Computer Science, 1993, 702:~425--439.
\bibitem{Zhuk_Le_Bars} M.E. Zhukovskii, {\it Logical laws for short existential monadic second order sentences about graphs}, arXiv:1712.06168, 2017.
\end{thebibliography}
\end{document}